\documentclass{llncs}
\usepackage{bbold,amssymb}
\usepackage{graphicx}
\spnewtheorem{pty}{Property}{\bfseries}{\itshape}

\newcommand{\fig}[3]{\begin{figure}[!hb]%
\begin{center}%
\includegraphics[scale=#3]{#1}%
\caption{#2}\label{#1}%
\end{center}%
\end{figure}%
}

\newcommand{\ignore}[1]{}

\let\kl\mathcal
\newcommand{\N}{\mathbb N}
\newcommand{\Z}{\mathbb Z}

\let\impl\Rightarrow
\let\equi\Leftrightarrow
\newcommand{\az}{A^\Z}
\newcommand{\an}{A^\N}
\newcommand{\lang}{\kl L}
\newcommand{\orb}{\kl{O}}

\newcommand{\sett}[2]{\left\{\left.#1\vphantom{#2}\right|#2\,\right\}}
\newcommand{\set}[3]{\sett{#1\in#2}{#3}}
\newcommand{\interv}[2]{\left[#1,#2\right]}
\newcommand{\ointerv}[2]{\interv{#1}{#2-1}}
\newcommand{\oointerv}[2]{\interv{#1+1}{#2-1}}
\newcommand{\subi}[2]{_{\interv{#1}{#2}}}
\newcommand{\suboi}[2]{_{\ointerv{#1}{#2}}}
\newcommand{\subinf}[1]{_{\left[0,+\infty\right)}}

\newcommand{\eqn}[1]{\begin{displaymath}\everymath{\displaystyle\everymath{}}\begin{array}{rclcl}#1\end{array}\end{displaymath}}
\newcommand{\both}[1]{\left\{\everymath{\displaystyle\everymath{}}\begin{array}{l}#1\end{array}\right.}
\newcommand{\soit}[1]{\left|\everymath{\displaystyle\everymath{}}\begin{array}{ll}#1\end{array}\right.}
\newcommand{\appl}[5]{#1:\begin{array}{rcl}#2&\to&#3\\#4&\mapsto&#5\end{array}}
\newcommand{\ou}{\textrm{ where }}
 \newcommand{\et}{\textrm{ and }}
 
 \newcommand{\si}{\textrm{ if }}
 \newcommand{\sinon}{\textrm{ otherwise}}
\newcommand{\mm}[1]{\textrm{ #1 }}

\newcommand{\length}[1]{\left|#1\right|}
\let\card\length
\let\abs\card

\newenvironment{prf}{\begin{proof}}{~\qed\end{proof}}

\newcommand{\xp}[3]{\vphantom{#1}^{#2}_{#3}#1}
\newcommand{\compl}[1]{#1^C}
\newcommand{\ie}{i.e.\ }
\newcommand{\defeq}{=}

\newcommand{\fact}{\sqsubset}

\newcommand{\eqref}[1]{(\ref{#1})}

\title{Sofic trace subshift of a cellular automaton%
\thanks{This work has been supported
by the Interlink/MIUR project ``Cellular Automata: Topological
Properties, Chaos and Associated Formal Languages", by the ANR
Blanc ``Projet Sycomore''.}
}
\author{Julien Cervelle\inst1\and Enrico Formenti\inst2\thanks{Corresponding author.}\and Pierre Guillon\inst1}
\institute{Institut Gaspard Monge, Universit\'e de Marne la Vall\'ee\\
77454 Marne la Vall\'ee Cedex 2, France.\\
\texttt{$\{$julien.cervelle,pierre.guillon$\}$@univ-mlv.fr}\\
\and Laboratoire I3S, Universit\'e de Nice-Sophia Antipolis\\
2000 Route des Lucioles, 06903 Sophia Antipolis, France.\\
\texttt{enrico.formenti@unice.fr}
}
\begin{document}

\maketitle

\begin{abstract}
The trace subshift of a cellular automaton is the subshift of all possible columns that may appear in a space-time diagram. 
In this paper we study conditions for a sofic subshift to be the trace of a cellular automaton. 
\end{abstract}


\section{Introduction}
Cellular automata are well-known formal models for 
complex systems.
They are used in a huge variety of different scientific fields including mathematics, physics and computer science.

A \emph{cellular automaton} (CA) consists in an infinite number
of identical \emph{cells} arranged on a regular lattice indexed by $\Z$.
Each cell is a finite automaton which \emph{state} takes value in
a finite set $A$. All cells evolve synchronously according to their own
state and those of their neighbors.

The study and classification of the evolutions of
cellular automata is one of the standing open problems in the field
\cite{gilman88,hurley90,dubacq99,durand03,cyu,catta}.
Indeed, the simple definition of CA contrasts the wide variety of 
their evolutions. An interesting idea is to classify these behaviors
according to some notion of complexity. Of course, the word
``complexity'' means different things to different researchers.
For this reason, in literature one finds classifications according
to topological entropy, measure theory, dimension theory, 
attractors, algorithmic complexity, etc.

In this paper we follow a formal languages approach.
Each CA is associated with a language. The idea is that the more
complex is the language, the more complex is the automaton.

The associated language is defined as follows
(see Section~\ref{sec:defs} for more precise definitions). 
Each CA can be seen as a discrete dynamical system $(\az,F)$, where $F$ is the global function. Let $\beta\defeq\beta_1,\beta_2,\ldots,\beta_k$ a finite partition of $\az$. Then an orbit of initial condition $x$, $\orb_F(x)=(x,F(x),\ldots,F^n(x),\ldots)$ can be associated with the infinite word $w$ such that
$\forall n\in\N,w_n=i$ if $F^n(x)\in\beta_i$. Then, $w$ is the
\emph{$\beta$-trace} of $F$ with initial condition $x$. The \emph{$\beta$-trace set} $\Sigma$ of $F$ is the set of $\beta$-traces with all possible
initial conditions. Remark that when $\az$ is endowed with
the Cantor topology (see Section~\ref{sec:defs}), $\Sigma$ is a 
closed shift-stable set \ie a subshift. The language $\lang(\Sigma)$ of factors occurring in configurations of $\Sigma$ is the language associated with $(\az,F)$. 

In~\cite{kurka97}, K\r{u}rka classified factor subshifts of CA 
according to their language complexity. He devised three classes:
bounded periodic; regular but not bounded periodic; not regular. 
In this paper we address a somewhat complementary question,
namely, given a subshift $\Sigma$ of a certain language complexity
we wonder if it can be the trace of a CA.

Motivations come both from classical symbolic dynamics but also
from physics. Indeed, when observing natural phenomena due to
physical constraints, one can keep trace only of a finite number of measurements. 
This set of measurements, usually, takes into account only a 
minor part of the parameters ruling the phenomenon under
investigation.
Hence, to some extent, what is observed is the ``trace'' 
of the phenomenon left on the ``instruments'' rather than the whole phenomenon in its globality.

It is a very important issue (Galilean
principle) to find a formal model which can reproduce the observed
trace. Restating those reasonings in our context: given a subshift
$\Sigma$, one wonders which discrete dynamical system can 
produce it. In particular, one can ask if there exists a CA having
$\Sigma$ as a trace.

Giving a complete answer to this question seems very hard. In this
paper we give some sufficient conditions for a regular language
to be traceable. The proof is constructive. We believe that the
construction of the CA is of some interest in its own.    
\medskip

Because of the lack of space some of the proofs are omitted. 
They can be found in the Appendices.

\section{Definitions}\label{sec:defs}
Let $\N^*=\N\setminus\{0\}$. For $i,j\in\N$ with $i\le j$,
$\interv ij$ denotes the set of integers between $i$ and $j$.  For
any function $F$ from $\az$ into itself, $F^n$ denotes the $n$-fold
composition of $F$ with itself. A set $S\subset\az$ is
\emph{$F$-stable} if $F(S)\subset S$.

\paragraph{Languages.}
Let $A$ be a finite \emph{alphabet} with at least two \emph{letters}.
A \emph{word} is a finite sequence of letters $w\defeq w_0\ldots
w_{\length w-1}\in A^*$. Its \emph{reverse} $\bar w$ is $w_{\length
w-1}\ldots w_0$ and its \emph{rotation} $\gamma(w)$ is $w_1\ldots
w_{\length w-1}w_0$.  A \emph{factor} of a word $w\defeq w_0\ldots
w_{\length w-1}\in A^*$ is a word $w\subi ij\defeq w_i\ldots w_j$, for
$0\le i\le j<\length w$. We note $w\subi ij\fact w$.  The empty word
is denoted by $\varepsilon$.  Given two languages $C,D\subset A^*$,
$CD$ denotes their concatenation, $C+D$ their union, $\displaystyle C^*\defeq\bigcup_{n\in\N}C^n$, and $C^\omega=\set z\an{\forall
j\in\N,\exists k\ge j,z\suboi0k\in C^*}$.  When no confusion is
possible, given a word $w$, we also denote $w$ the language $\{w\}$.

\paragraph{Configurations.}
A \emph{configuration} is a biinfinite sequence of letters $x\in\az$. The set $\az$ of configurations is the \emph{phase space}.
The definition of \emph{factor} can be naturally extended to configurations: for $x\in\az$ and $i\le j$, 
$x\subi ij\defeq x_i\ldots x_j\fact x$.
If $u\in A^*$, then $u^\omega$ is the infinite word consisting in periodic repetitions of $u$, and $\xp u\omega{}^\omega$ is the configuration consisting in periodic repetitions of $u$.

\paragraph{Topology.}
We endow the phase space with the \emph{Cantor topology}.
A base for open sets is given by cylinders. 
For $j,k\in\N$ and a finite set $W$ of words of length $j$, we will 
note $[W]_k$ the \emph{cylinder} $\set w\az{w\suboi k{k+j}\in W}$.
$\compl{[W]_k}$ is the complement of the cylinder $[W]_k$.

\paragraph{Cellular automata.}
A (one-dimensional) \emph{cellular automaton} is a parallel
synchronous computation model consisting in cells distributed over
a regular lattice $\Z$. 
Each cell has a state in the finite alphabet $A$, which evolves depending on the state of their neighbors according to a \emph{local rule} $f:A^d\to A$, where $m\in\Z$ and $d\in\N^*$ are the \emph{anchor} and the \emph{diameter} of the CA, respectively.
The \emph{global function} of the CA is $F:\az\to\az$ such that $F(x)_i=f(x\subi{i-m}{i-m+d})$ for every $x\in\az$ and $i\in\Z$.
The \emph{space-time diagram} of initial configuration $x\in\az$ is the sequence of the configurations of the orbit $(F^j(x))_{j\in\N}$.
Usually they are graphically represented by a two-dimensional
diagram like in Figure~\ref{trace}.

The \emph{shift map} $\sigma:\az\to\az$ is a particular CA global function defined by $\sigma(x)_i=x_{i+1}$ for every $x\in\az$ and $i\in\Z$.
According to the Hedlund theorem~\cite{hedlund69}, the global functions of CA are exactly the continuous
self-maps of $\az$ commuting with the shift map.

Any local rule $f$ of a CA can be extended naturally to an application
on words $f(w)=(f(w\suboi i{i+d}))_{0\le i\le\length w-d}$, for all $w\in
A^*A^d$.

\paragraph{Dynamical systems.}
A \emph{dynamical system} is a couple $(X,F)$ where $X$ is a set called the \emph{phase space}, and $F:X\to X$ is a continuous self-map.
$(Y,F)$ is a \emph{subsystem} of $(X,F)$ if $Y$ is a closed 
$F$-stable subset of $X$. A set $Y$ is \emph{$F$-stable} 
if $F(Y)\subset Y$.

\paragraph{Morphisms.}
A \emph{morphism} of the dynamical system $(X,F)$ into the dynamical system $(Y,G)$ is a continuous map $\phi:X\to Y$ such 
that $\phi\circ F=G\circ\phi$. A \emph{conjugacy} (resp. a 
\emph{factorization}) is a bijective (resp. surjective) morphism; in that
case we say that $(Y,G)$ is conjugate to (resp. a factor of) $(X,F)$ 
(we expect the reader not to confuse between the use of the word
``factor'' in the dynamical system context and the language theory one).

\ignore{
The following remark will be very useful in the sequel.

\begin{remark}\label{subconj}
If $(X,F)$ is a factor of (conjugate to) $(Y,G)$, then any subsystem
$(X',F)$ is a factor of (conjugate to)
$(\phi(X'),G)$, where $\phi$ is the factorization (conjugacy) map and
$X'\subset X$.
\end{remark}
}

\paragraph{Subshifts.}
The \emph{onesided shift}, also noted $\sigma$ is the self-map of
$\an$ such that $\sigma(z)_i=z_{i+1}$, for every $z\in\an$ and
$i\in\N$.  A \emph{onesided subshift} (or simply a \emph{subshift}) 
$\Sigma\subset A^\N$ is a $\sigma$-stable closed set
of infinite words. The language of $\Sigma$ is $\lang(\Sigma)\defeq\set w{A^*}{\exists z\in\Sigma,w\fact z}$, and characterizes $\Sigma$, since $\Sigma=\set z{A^\N}{\forall w\fact z,w\in\lang(\Sigma)}$.
The \emph{alphabet} of the subshift $\Sigma$ is the set 
$\set aA{\exists z\in\Sigma,az\in\Sigma}$ \ie the set of letters that
appear in infinite words belonging to $\Sigma$.
A subshift $\Sigma$ is \emph{transitive} if for every words $u,v\in\lang(\Sigma)$, there exists $w\in A^*$ such that $uwv\in\lang(\Sigma)$.
A subshift can also be characterized by a language 
$\mathcal{F}\subset A^*$ of \emph{forbidden words}, 
\ie such that $\Sigma=\set z\an{\forall u\in\mathcal{F},x\not\fact z}$. 
A subshift is of \emph{finite type} (SFT for short) if it has a finite
language of forbidden words. It is a $k$-SFT (for $k\in\N$) if it has 
a finite set of forbidden words of length $k$.
A subshift $\Sigma$ is \emph{sofic} if $\lang(\Sigma)$ is a
regular language.
The following characterization of sofic subshifts will be very useful in the sequel.
For more about subshifts, see for instance~\cite{lind95}.
\begin{theorem}[Weiss \cite{weiss73}]\label{factsft}
  A subshift is sofic if and only if it is a factor of a SFT. 
\end{theorem}
Hedlund's theorem can be extended as follows.
\begin{theorem}[Hedlund \cite{hedlund69}]
   A function $\phi$ is a morphism of a subshift $(X,\sigma)$ on 
  alphabet $A$ into a subshift
  $(Y,\sigma)$ on alphabet $B$ if and only if there is a radius
  $r\in\N$ and a local rule $f:A^{r+1}\to B$ such that $\forall
  j\in\N,\forall x\in X,\phi(x)_j=f(x_j\ldots x_{j+r})$ (we say $\phi$
  is an $r$-block map).
\end{theorem}

\ignore{
The following easy result will be used in the sequel.

\begin{proposition}\label{chalpha}
If $A$ and $B$ are two alphabets such that $\card B\ge\card A$, then any subshift on $A$ is conjugate to some subshift on $B$.
\end{proposition}
\begin{prf}
There is an injection $f$ from $A$ into $B$, which is a 1-block local rule. If $\Sigma$ is a subshift on $A$, then it is conjugate to subshift $\sett{(f(z_j))_{j\in\N}}{z\in\Sigma}$ on $B$.
\end{prf}
}


%
\section{Traces}
In this section we define the main notion introduced in the paper,
namely, the trace of a CA and the traceability of a subshift.
Moreover, we give a simple necessary condition for a subshift
being traceable.

\fig{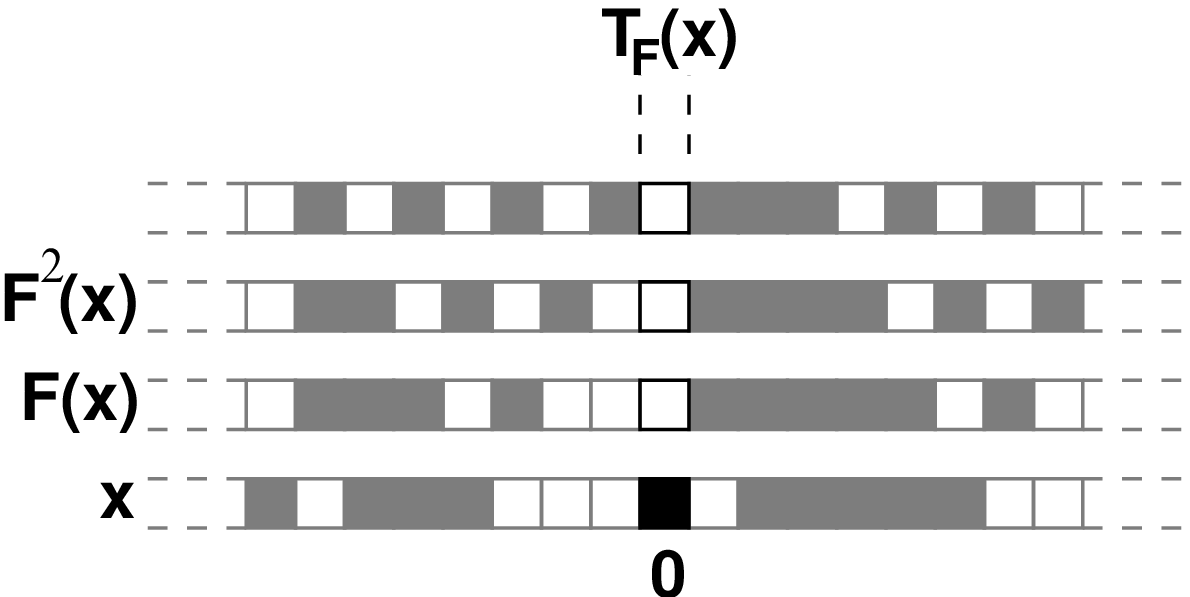}{The trace seen on the space-time diagram.}{.4}
\begin{definition}[Trace]
Given a CA $F$, the \emph{trace applications} are defined for $k\in\Z$ by $T_F^k(x)\defeq(F^j(x))_{j\in\N}$. In other words, $T_F^k(x)$ is the $k^{\rm th}$ column of the space-time diagram of initial configuration $x$ (see Figure~\ref{trace}). We note $T_F=T^0_F$. We say that $T_F(x)$ is the trace of $F$ with initial condition $x$.
\end{definition}
The study of trace applications can be reduced to the study of $T_F$, because of shift-invariance of CA. 

It can be noticed that this notion of trace corresponds to the $\beta$-trace as defined in the introduction, with $\beta\defeq\sett{[a]}{a\in A}$ being the partition of $\az$ into cylinders of width $1$.

\begin{definition}[Traceability]
The \emph{trace subshift} of a CA $F$ is $\tau(F)\defeq T_F(\az)$.
It is a factor subshift of $(\az,F)$, since $T_F$ is continuous and commutes with $\sigma$. A subshift $\Sigma$ is \emph{traceable} if there exists a CA $F$ for which $\Sigma=\tau(F)$.
\end{definition}

We begin with a condition for the traceability of a subshift.
Proposition~\ref{t0} proves that it is necessary.
\begin{definition}[T0 subshift]
A subshift is \emph{T0} if it includes a $2$-SFT with the same alphabet.
\end{definition}

\begin{proposition}\label{eqt0}
A subshift is T0 if and only if there exists a map $\phi:A\to A$ such 
that for every letter $a\in A$, $(\phi^j(a))_{j\in\N}\in\Sigma$.
\end{proposition}
\begin{prf}
If $\forall a\in A,(\phi^j(a))_{j\in\N}\in\Gamma$, then $\sett{(\phi^j(a))_{j\in\N}}{a\in A}$ is a $2$-SFT of alphabet $A$ and set of forbidden words $\bigcup_{a\in A} a(A\setminus\{\phi(a)\})$.
Conversely, consider a  $2$-SFT $\Gamma$ of alphabet $A$. Define $\phi(a)=b$ such that $ab\in\lang(\Gamma)$. Since $\Gamma$ is a $2$-SFT, we have that $\forall a\in A,(\phi^j(a))_{j\in\N}\in\Gamma$.
\end{prf}

\begin{example}
Consider the following subshifts :
\begin{itemize}
\item $\Sigma=(1+\varepsilon)(01)^\omega$; it is finite T0 (with $\phi(0)=1$ and $\phi(1)=0$); 
\item $\Sigma=0^\omega+0^*1^\omega$; it is infinite T0 (with $\phi(0)=\phi(1)=1$);
\item $\Sigma=(01+1+\varepsilon)(001)^\omega$; it is finite 
but not T0.
\end{itemize}\end{example}

\begin{proposition}\label{t0}
The trace subshift of a CA is T0.
\end{proposition}
\begin{prf}
Consider a CA  $F$. For any $a\in A$, define $\phi(a)$ as 
$F(\xp a\omega{}^\omega)_0$. 
Then $(\phi^j(a))_{j\in\N}=T_F(\xp a\omega{}^\omega)\in\tau(F)$. 
By Proposition~\ref{eqt0}, $\tau(F)$ is T0.
\end{prf}

\begin{theorem}\label{sft2}
  Any $2$-SFT is traceable.
\end{theorem}
\begin{prf}
  Consider a 2-SFT $\Sigma$ and let $A$ be its alphabet.
  For each $a\in A$, define $\phi(a)=b$
  such that $ab\in\lang(\Sigma)$. 
  Consider the CA of anchor $0$, diameter $2$ and local rule $f$ defined by $f(x_0,x_1)=x_1$ if $x_0x_1\in\lang(\Sigma)$, and $\phi(x_0)$ otherwise.
If $x\in\az$, then the definition of the rule gives that every factor of length $2$ of its trace is in $\lang(\Sigma)$. 
Conversely, if $z\in\Sigma$, then we can see by induction that $\Sigma$ is the trace of $F$ with initial condition $x$ as soon as $x\subinf0=z$.
\end{prf}
\section{$k$-traceability}\label{ktrace}
In order to establish finer results we first need a weaker 
condition for traceability, namely $k$-traceability.
A subshift of alphabet $A$ is $k$-traceable if it is the set of
columns of a CA on the alphabet $A^k$. The difference between
a CA tracing a subshift of alphabet $A$ and one tracing a subshift
of alphabet $A^k$ is that the rule of the latter can use the
knowledge of the position of a letter of $A$ in a word over
$A^k$. This results in a much simpler construction.
\smallskip

\paragraph{Notation.}
If $\Sigma$ is a subshift on an alphabet $B\subset A^k$, and $q\in\ointerv0k$, then the $q^\mathrm{th}$ projection is defined as
\[\appl{\pi_q}{B^\N}\an{(z_j)_{j\in\N}}{((z_j)_q)_{j\in\N}}\enspace.\]
We also note $\pi(\Sigma)\defeq\bigcup_{0\le q<k}\pi_q(\Sigma)$, which is a subshift on $A$.

In this section we will limit our study to onesided CA, \ie with anchor $0$.
\begin{definition}[$k$-traceability]
Given a CA $F$ on the alphabet $B\subset A^k$, the 
\emph{$k$-trace subshift} is defined by $\stackrel\circ\tau(F)\defeq\displaystyle\bigcup_{0\le q<k}\sett{((F^j(x)_0)_q)_{j\in\N}}{x\in B^\Z}=\pi(\tau(F))$. 
A subshift is \emph{$k$-traceable} if it is the $k$-trace of a onesided CA on the alphabet $B\subset A^k$.
\end{definition}

Similarly to what done in the previous section we give a necessary
condition for being $k$-traceable. 

\begin{definition}[T1 subshift]
A subshift $\Sigma$ on the alphabet $A$ is \emph{T1} if there exists $k\in\N^*$ and a $2$-SFT $\Gamma\subset(A^k)^\N$, such that $\pi_0(\Gamma)=\pi(\Gamma)=\Sigma$ (in particular, $\Sigma$ is a factor of $\Gamma$).
\end{definition}

\begin{theorem}\label{t1}
A T1 subshift is $k$-traceable for some $k\in\N^*$.
\end{theorem}
\begin{prf}
By Theorem~\ref{sft2}, the corresponding 
$\Gamma\subset(A^k)^\N$ is
the trace of a CA $F$ on some alphabet $B\subset A^k$. 
Hence, $\stackrel\circ\tau(F)=\pi(\tau(F))=\pi(\Gamma)=\Sigma$.
\end{prf}

\begin{example}
Consider the subshift $\Sigma=(01+1+\varepsilon)(001)^\omega$. It is T1 (define the $2$-SFT $\Gamma=(uvw)^{\omega}$ on the alphabet $B=\{u,v,w\}\subset A^3$, where $u=001$, $v=010$ and $w=100$). It is thus $3$-traceable, but not traceable since it is not T0.
\end{example}

\ignore{
\begin{remark}
  Any traceable subshift is the $k$-trace of a CA with $m=1$ 
  and $d=3$.
\end{remark}
\begin{prf}
  Consider a CA $F$ on alphabet $A$ of anchor $m$, diameter $d$ 
  and local rule $f:A^d\to A$.
If $k\defeq\max(1,m,d-m-1)$, we can transform it into a CA $\tilde F$ with $m=1$, $d=3$ and local rule
\[\appl{\tilde f}{(A^k)^3}A{(u_0\ldots u_{k-1},u_k\ldots u_{2k-1},u_{2k}\ldots
  u_{3k-1})}{(f(u_{q-m},\ldots u_{q-m+d-1}))_{k\le q<2k}}\]
Then $\stackrel\circ\tau(\tilde F)=\tau(F)$.
\end{prf}
}

\begin{theorem}\label{thsftt1}
  Any SFT is T1.
\end{theorem}
\begin{prf}
Let $\Sigma$ be a $k$-SFT for some $k\in\N^*$. 
$\Gamma\defeq\sett{(z\subi j{j+k-1})_{j\in\N}}{z\in\Sigma}$ is a $2$-SFT, and $\pi(\Gamma)=\bigcup_{0\le q<k}\sett{(z_{j+q})_{j\in\N}}{z\in\Sigma}
=\bigcup_{0\le q<k}\sigma^q(\Sigma)
=\Sigma=\pi_0(\Gamma)$.
\end{prf}

\ignore{\fig{sft0}{CA $3$-tracing the $3$-SFT of forbidden language $\{111\}$.}{.4}}

This result allows us to prove the next proposition, which is a less restrictive condition for being T1.
\begin{proposition}\label{factt1}
A subshift $\Sigma$ is T1 if and only if it is a factor of a SFT 
$\Gamma$ on alphabet $A^k$ for some $k\in\N^*$ such that
$\pi(\Gamma)\subset\Sigma$.
\end{proposition}

Now we extend the results on $k$-traceability to sofic subshifts (with some additional properties).

\begin{definition}[T2 subshift]
A subshift is \emph{T2} if it is sofic and includes an infinite transitive subshift.
\end{definition}

\begin{theorem}\label{tht3}
Any T2 subshift is T1.
\end{theorem}
The proof of Theorem~\ref{tht3} is given using the following lemmata.

\begin{lemma}\label{choix}
A sofic transitive subshift on alphabet $A$ is infinite if and only if
for all $n\ge2$, it includes a subshift $B^\omega$ with $B\subset
A^k$, $\card B\ge n$ and some $k\in\N^*$.
\end{lemma}

\begin{lemma}\label{artificiel}
If $\Sigma$ is a factor subshift of a SFT $\Gamma$ on alphabet $B$
such that $B^\omega\subset\Sigma$, then $\Sigma$ is T1.
\end{lemma}

\begin{prf}[of Theorem~\ref{tht3}]
Let $\Sigma$ a T2 subshift. From Theorem~\ref{factsft}, it is a factor
of a SFT $\Gamma$. Thanks to Lemma~\ref{choix}, there is an arbitrarily large set of words $B$ on $A$ such that $B^\omega\subset\Sigma$, so we can assume without loss of generality that $\Gamma$ is a subshift on such an alphabet $B\subset A^k$ for some $k\in\N^*$. From Lemma~\ref{artificiel}, we conclude that $\Sigma$ is T1.
\end{prf}
\section{From $k$-trace to trace}\label{k2trace}
%
In the previous section, we gave a sufficient condition for a
particular subshift $\Sigma$ to be $k$-traced by a CA $G$ on an alphabet $B\subset A^k$.
In this section, we show how to simulate $G$ with another CA, on
alphabet $A$, in such way that its trace is $\Sigma$. This can be done if we add a further condition to our subshift.

\begin{definition}[T3 subshift]\label{t2def}
A subshift $\Sigma$ is \emph{T3} if there is a map $\phi:A\to A$ such that for every letter $a\in A,(\phi^j(a))_{j\in\N}\in\Sigma$ (it is T0) and there is a word $w\in A^*\setminus\phi(A)^*$ such that $w^\omega\in\Sigma$.
\end{definition}
\begin{example}
Consider the following subshifts.
\begin{itemize}
\item$\Sigma\defeq(1+\varepsilon)(01)^\omega$ is not T3.
\item$\Sigma'\defeq0^\omega+0^*1^\omega$ is T3 (with $\phi(0)=\phi(1)=1$ and $w=0$).
\end{itemize}\end{example}

\begin{theorem}\label{t2}
Any T3 $k$-traceable subshift (for some $k\in\N^*$) is traceable.
\end{theorem}
This section presents a sketch of the proof of Theorem~\ref{t2}.
Remark that it is well known that a CA on any alphabet can be simulated by a CA on any other alphabet (with at least two letters), provided that its diameter is wide enough. In particular, any CA on $B\subset A^k$ can be simulated by a CA on $A$. Each cell can see its neighborhood as words of $A^k$ and evolve accordingly. The problem is that all cells must have the same local rule, so they have to find from the neighborhood which \emph{column} of the $A^k$ simulation they are representing.
This is usually done using a special \emph{border} word to delimit the words of $A^k$.

In this section, $\Sigma$ denotes a T3 subshift on alphabet which is $k$-traceable by a onesided CA $G$. Let $\phi$ and $w$ be as in Definition~\ref{t2def}. Assume $G$ has diameter $2$ (the construction can easily be generalized) and local rule $g:B^2\to B$.

In order to achieve the simulation, we first define border words to
delimit $A^k$ cells. We have two \emph{execution modes}: a
\emph{simulation mode} will simulate properly the execution of the CA
on alphabet $B$, and a \emph{default mode} will be applied if the
neighborhood contains invalid information.
This adds some issues: 
default evolution must be in $\Sigma$; border evolution must also
evolve according to $\Sigma$; and we have to ensure that when a mode
is applied to a cell, the same mode keeps being applied there in the
following generations, since a change of mode would produce an invalid
trace. These problems will be solved in the three following subsections.

\subsection{Borders}
In order to make our simulations, we need to delimit computation
zones. This is obtained by using some special words called borders and
defined as follows:
\[\Upsilon\defeq\sett{a^{\length w}v\overline va^{k+3\length
    w}}{a\in\phi(A),v\in\orb_\gamma(w)}\subset A^l,\]
where $\orb_\gamma(w)\defeq\sett{\gamma^q(w)}{0\le q<\length w}$, and $l\defeq k+6\length w$.

Borders have the property that they cannot have a too wide overlap. 

\paragraph{Border evolution.}
The border words of $\Upsilon$ must have an evolution in $\Sigma$. The following rule (of diameter $1$) respects that condition:
\[\appl{\Delta_\Upsilon}\Upsilon\Upsilon{a^{\length w}v\overline va^{k+3\length w}}{\phi(a)^{\length w}\gamma(v)\overline{\gamma(v)}\phi(a)^{k+3\length w}}\enspace.\]
\ignore{\fig{border}{a border element for $\length w=4$, and two iterations of $\Delta_\Upsilon$.}{.4}}

\paragraph{Macrocells.}
We will decompose our configurations into \emph{macrocells}. A macrocell
is the concatenation of a border word and a valid word of $B$. We can
simulate the local rule $g$ by a \emph{macroevolution rule}
(local rule on macrocells of $B\Upsilon\subset A^h$, where $h\defeq
k+l=2k+6\length w$, and of diameter $2$):
\[\appl\Delta{(B\Upsilon)^2}{B\Upsilon}{(u,v)}{g(u\suboi0k,v\suboi0k)\Delta_\Upsilon(u\suboi kh)}\enspace.\]

\subsection{Default mode}
Our CA will work as follows. Valid zones (with macrocells), which evolve according to the macroevolution rule $\Delta$ so that they remain valid zones. Invalid zones run a microdefault mode so that they remain invalid zones. Nevertheless, frontiers between the zones must not move. In the frontiers, a macrodefault mode is applied in order for a macrocell to have the opportunity to evolve without taking into account its neighbors; that way, each cell will keep the same execution mode.

\paragraph{Macrodefault mode.}
To do so, we extend the macroevolution to a function on $\Theta A^h$, where $\Theta\defeq B\Upsilon A^h\setminus\displaystyle\bigcup_{0<i<h}A^iB\Upsilon A^{h-i}$ because it does not take into account overlapping macrocells. This is crucial in order to define a local rule. If the central macrocell has a neighbor macrocell in $\Theta$, we apply a simulation step of the CA. Otherwise, we evolve as a macrodefault mode (simulation from a monochromatic configuration):
\[\appl\Delta{\Theta A^h}{B\Upsilon}u{\soit{g(u\suboi0k,u\suboi h{h+k})\Delta_\Upsilon(u\suboi kh)\si u\in B\Upsilon\Theta\\g(u\suboi0k,u\suboi0k)\Delta_\Upsilon(u\suboi kh)\sinon}}\enspace.\]

\ignore{\fig{sft}{CA tracing the $3$-SFT of forbidden language $\{111\}$.}{.4}}

\paragraph{Microdefault mode.}
The function $\phi$ (corresponding to the fact $\Sigma$ is T0) allows to define a microdefault mode for a neighborhood that does not contain any macrocell. We are now able to transform the function $\Delta$ into a local rule on $A$.
Indeed, we can define, for anchor $m\defeq h-1$ and diameter $d\defeq 3h-1$:
\[\appl f{A^d}Aw{\soit{\Delta(u)_i&\si w\in A^{m-i}uA^i,\ou u\in\Theta A^h,i\in\ointerv0h\\\phi(w_0)&\sinon}}\]
since such an integer $i$, and such a word $u$ would be unique (from the construction of $\Theta$).
This local rule is such that $f(A^muA^m)=\Delta(u)$ for every $u\in\Theta A^h$, which is what we wanted: it can simulate in one step the behavior of our CA on $B$. Let $F$ be the corresponding global rule.

The following lemma guarantees that no column changes its evolution mode. 
\begin{lemma}\label{stabl}
The preimage of cylinder $[B\Upsilon]$ is cylinder $[\Theta]$. Moreover, cylinder $[\Theta]$ and its complementary $\compl{[\Theta]}$ are $F$-stable (in particular, we cannot create a border).
\end{lemma}

A configuration which is a valid encoding of some $y\in B^\Z$ (simulation mode), then its trace is some projection of the trace of $y$. Otherwise, microdefault and microdefault mode also produce a trace which is in $\Sigma$.
This concludes the proof of Theorem~\ref{t2}.

\section{Examples}

\begin{example}[Finite untraceable T1 subshift]\label{ctrex}
No CA traces subshift $\Sigma\defeq\{0^\omega,(01)^\omega,(10)^\omega\}$, even though it is T1.
\end{example}


\begin{example}[T1, T3, non-SFT, non-T2 subshift]\label{T1T3nonT2}
The subshift $\Sigma=(0^*1+1^*)0^\omega$ is neither a SFT nor T2, but it is T1 and T3. Hence, by Theorems~\ref{t1} and \ref{t2} it is traceable.
\end{example}

\begin{example}[Traceable non-T3 subshift]
Let $f$ be the local rule of anchor $3$ and diameter $7$ such that $f(u_{-3}000111)=1$, $f(000111u_3)=0$, $f(u_{-3}001011)=0$, $f(001011u_3)=1$, and $f(u)=u_0$ otherwise.
The {trace subshift} of the corresponding CA is $\tau(F)=\{0^\omega,(10)^\omega,(01)^\omega,1^\omega\}$. In this case, $\tau(F)$ is finite but not T3.
\end{example}

\begin{example}[Traceable non-sofic subshift]\label{ctrex2}
Let $F$ be the CA on alphabet $A\defeq\{b,r,l,w\}$ (the white, the
right, the left and the wall particles, respectively) defined by the
following local rule $f$ of anchor $1$ and diameter $3$:
\[\begin{array}{|r|c|c|c|c|c|c|c|c|c|}\hline
x_{-1}x_0x_1&rl?&?rl&r?l&?w?&?rw&wl?&r??&??l&???\\
\hline
f(x_{-1}x_0x_1)&w&w&w&w&l&r&r&l&b\\
\hline\end{array}\]
where $?$ stands for any letter in $A$ and the first applicable rule is used (left to right).
Then, $\tau(F)$ is not sofic.
\end{example}

\section{Putting things together}
In this paper, we have given sufficient conditions for a subshift to be the trace of a CA.
\begin{eqnarray*}
\mm{T1}&\impl&\mm{$k$-traceable (Theorem~\ref{t1})}\\
\mm{SFT}&\impl&\mm{T1 (Theorem~\ref{thsftt1})}\\
\mm{T2}&\impl&\mm{T1 (Theorem~\ref{tht3})}\\
\mm{T3 and $k$-traceable}&\impl&\mm{traceable (Theorem~\ref{t2})}
\end{eqnarray*}
The following summarizes all these results:
\begin{theorem}
Any T3+T1, T3 SFT, or T3+T2 subshift is traceable.
\end{theorem}

The present result follows other works on the structure that the trace of a CA can have. Here, we take the problem the other way around: we construct a CA that traces a particular kind of subshifts. Though we do not have a necessary and sufficient condition, Conditions T0, T1, T2 and T3 are a first step toward a better understanding of what makes a sofic subshift traceable or not.

Moreover, we expect our construction to be generalizable to weaker conditions. Nevertheless, the non-sofic case is still obscure. The general feeling is that it needs a completely different approach.

In \cite{kurka97}, it is proved that every factor subshift of a CA $F$ is a factor of some $\beta$-trace, where $\beta$ is a partition of $\az$, and every $\beta$-trace is a factor of some \emph{column factor} (\ie a subshift $\sett{(F^j(x)\suboi0q)_{j\in\N}}{x\in\az}$, for some $q\in\N$). Hence we can wonder now whether this kind of result can be generalized, in particular to the \emph{canonical factor} $\sett{(F^j(x)\suboi0d)_{j\in\N}}{x\in\az}$ of the CA.

\bibliographystyle{splncs}
\bibliography{tra}

\begin{thebibliography}{10}

\bibitem{gilman88}
Gilman, R.H.:
\newblock Classes of linear automata.
\newblock Erg. Th. {\&} Dyn. Sys. \textbf{7} (1988)  105--118

\bibitem{hurley90}
Hurley, M.:
\newblock Attractors in cellular automata.
\newblock Erg. Th. {\&} Dyn. Sys. \textbf{10} (1990)  131--140

\bibitem{dubacq99}
Dubacq, J.C., Durand, B., Formenti, E.:
\newblock Kolmogorov complexity and cellular automata classification.
\newblock Th. Comp. Sci. \textbf{259}(1--2) (2001)  271--285

\bibitem{durand03}
Durand, B., Formenti, E., Varouchas, G.:
\newblock On undecidability of equicontinuity classification for cellular
  automata.
\newblock In Morvan, M., R{\'e}mila, E., eds.: DMCS'03. Volume~AB of DMTCS
  Proc., Disc. Math. and Th. Comp. Sci. (2003)  117--128

\bibitem{cyu}
Culik, K., Yu, S.:
\newblock Undecidability of cellular automata classification schemes.
\newblock Comp. Sys. \textbf{2} (1988)  177--190

\bibitem{catta}
Braga, G., Cattaneo, G., Flocchini, P., Vogliotti, C.Q.:
\newblock Pattern growth in elementary cellular automata.
\newblock Th. Comp. Sci. \textbf{145}(1--2) (1995)  1--26

\bibitem{kurka97}
K{\r u}rka, P.:
\newblock Languages, equicontinuity and attractors in cellular automata.
\newblock Erg. Th. {\&} Dyn. Sys. \textbf{17} (1997)  417--433

\bibitem{hedlund69}
Hedlund, G.A.:
\newblock Endomorphism and automorphism of the shift dynamical system.
\newblock Math. Sys. Theory \textbf{3} (1969)  320--375

\bibitem{lind95}
Marcus, B., Lind, D.:
\newblock An introduction to symbolic dynamics and coding.
\newblock Cambridge University Press (1995)

\bibitem{weiss73}
Weiss, B.:
\newblock Subshifts of finite type and sofic systems.
\newblock Monatshefte f{\"u}r Mathematik \textbf{77}(5) (1973)  462--474

\end{thebibliography}
\newpage
\appendix

\section{Proofs of Section~\ref{ktrace}}

\begin{prf}[of Proposition~\ref{factt1}]
Assume that $\Sigma$ is a T1 subshift; then there is a $2$-SFT 
$\Gamma$ such that $\pi_0(\Gamma)=\Sigma$. 
As $\pi_0$ is a factorization, $\Sigma$ is a factor of the SFT $\Gamma$.

For the converse implication, assume that $\Sigma$ is a subshift on alphabet $A$, $\Gamma$ a SFT on alphabet $B\subset A^k$, and $\phi:\Gamma\to\Sigma$ a factorization, where $k\in\N$ and $\pi(\Gamma)\subset\Sigma$. We can suppose without loss of generality that $\pi(\Gamma)=\pi_0(\Gamma)=\Sigma$, should we replace it by its conjugate $\set{(a_jw_j)_{j\in\N}}{(AB)^\N}{(w_j)_{j\in\N}\in\Gamma\et(a_j)_{j\in\N}=\phi((w_j)_{j\in\N})}$, which is still a SFT.
From Theorem~\ref{thsftt1}, $\Gamma$ is T1, \ie there exists a $2$-SFT $\Gamma'$ on alphabet $C\subset B^l$ for some $l\in\N$ such that $\pi_0(\Gamma')=\pi(\Gamma')=\Gamma$. We have just used the projections with respect to alphabet $B$ (recall that $C\subset B^l\subset(A^k)^l$), which projections with respect to $A$ itself are $\pi_0(\Gamma)=\pi(\Gamma)=\Sigma$. Hence, with respect to $A$, we get $\pi_0(\Gamma')=\pi(\Gamma')=\Sigma$, and $\Sigma$ is T1.
\end{prf}

\ignore{\begin{remark}
We can equivalently see a sofic subshift as the set of labels of infinite paths in a labeled directed graph. Moreover, a sofic subshift is transitive if the graph of the corresponding automaton is strongly connected.
\end{remark}}

\begin{prf}[of Lemma~\ref{choix}]
If $\card B\ge2$, then $B^\omega$ is infinite (since all words of $B$ have the same length). Conversely, if $\Sigma$ is an infinite transitive sofic subshift, and $(A,Q,\delta)$ a deterministic finite automaton which recognizes the regular language $\lang(\Sigma)$, then there is a state $q\in Q$ and two letters $a,b\in A$ such that $\delta(a,q)$ and $\delta(b,q)$ are defined (otherwise there would only be one possibility to extend the word from each state, and there would be at most as much infinite words in $\Sigma$ as initial states in $Q$). By transitivity, there is a word $\tilde u$ and a word $\tilde v$ such that $\delta(\tilde u,\delta(a,q))=\delta(\tilde v,\delta(b,q))=q$. If $n\defeq\length{a\tilde u}\cdot\length{b\tilde v}$, $u\defeq(a\tilde u)^{\length{b\tilde v}}$, $v\defeq(b\tilde v)^{\length{a\tilde u}}$, and $l$ is such that $2^l\ge n$, then $B=\{u,v\}^l\subset(A^n)^l$ is appropriate.
\end{prf}

\subsection{Proof of Lemma~\ref{artificiel}}
The function $H$ defined in the following lemma allows us to define a clock of period $2n$: the word $H(i)$ represents the $i$-th clock tick.
\begin{lemma}\label{clock}
Consider an alphabet $B\subset A^n$, with $\card B\ge2$ and $n\in\N^*$.
Then there is an injection $H:\ointerv0{2n}\to A^{3n}$ such that the subshift $\Gamma\defeq\sett{(H((q+j)\bmod2n))_{j\in\N}}{q\in\N}$
verifies $\pi(\Gamma)\subset B^\omega$.
\end{lemma}
\begin{prf}
\ignore{\fig{clock0}{The clock injection.}{.4}}
Let $u,v\in B$, suppose $u_0\ne v_0$ (other cases are obtained by rotation), and define:
\[\appl H{\ointerv0{2n}}{A^{3n}}h{\gamma^h(u)\gamma^h(uv)}\enspace.\]
Assume $H(h_1)=H(h_2)$ for some distinct $h_1,h_2\in\ointerv0n$. Let $h\defeq\min(h_1-h_2\bmod2n,h_2-h_1\bmod2n)\in\interv1n$. Then, by the definition of $H$, $\gamma^{h_1}(u)\gamma^{h_1}(uv)=\gamma^{h_2}(u)\gamma^{h_2}(uv)$. Hence, $\gamma^{h_1}(u)=\gamma^{h_2}(u)$ and $\gamma^{h_1}(uv)=\gamma^{h_2}(uv)$, giving $u=\gamma^h(u)$ and $uv=\gamma^h(uv)$.
Finally,
\[u_0=\gamma^h(u)_{n-h}=u_{n-h}=(uv)_{n-h}=\gamma^h(uv)_{n-h}=(uv)_n=v_0\enspace.\]
This is a contradiction, hence $H$ is injective.

Let $q\in\ointerv0{2n}$ and $p\in\ointerv0n$. Then,
\begin{eqnarray*}
\pi_p((H((q+j)\bmod2n))_{j\in\N})&=&(\gamma^{(q+j)\bmod2n}(u))_{j\in\N}
=(u_{(q+j)\bmod n})_{j\in\N}\\
&=&\sigma^q(u^\omega)\in B^\omega
\enspace.
\end{eqnarray*}
Similarly, taking the projection $p$ in $\ointerv n{3n}$ we find
\begin{eqnarray*}
\pi_p((H(q+j\bmod2n))_{j\in\N})&=&(\gamma^{q+j\bmod2n}(uv))_{j\in\N}
=((uv)_{q+j\bmod2n})_{j\in\N}\\
&=&\sigma^q((uv)^\omega)\in B^\omega
\enspace.\end{eqnarray*}
\end{prf}

\begin{lemma}\label{multi}
A full shift $(B^\N,\sigma)$, where $B\subset A^n$ and $n\in\N^*$, is
a factor of a SFT $\Psi$ on alphabet $A^{4n}$ such that $\pi(\Psi)\subset\orb_\sigma(B^\omega)$, seen as a subshift on $A$.
\end{lemma}
\begin{prf}
For $0\le q<n$, we define the set of \emph{time-$q$ encodings} as
\[Psi_q\defeq\set{((w_0)_j,\ldots,(w_{n-1})_j)_{j\in\N}}{(A^n)^\N}{\forall
  p\in\ointerv0n,w_p\in\sigma^{q-p\bmod n}(B^\omega)}\]
and the \emph{decoding function} as:
\[\appl{\phi_q}{\Psi_q}{B^\N}{((w_0)_j\ldots(w_{n-1})_j)_j\in\N}{((w_{q+j\bmod n})\suboi0n)_{j\in\N}}\enspace.\]
$\phi_q$ is surjective.
If we apply $\sigma$, we pass from a time-$q$ encoding to a time-$(q+1)$ encoding: $\sigma(\Psi_q)=\Psi_{q+1\bmod n}$.
In order to know which step we are in (which of the first $n$ columns to look at), we use a ``clock'' represented by the last $3n$ columns as described in Lemma~\ref{clock}.
We use the set:
\[\Psi\defeq\bigcup_{0\le q<2n}\set{(w_j,H(q+j\bmod2n))_{j\in\N}}{(A^{4n})^\N}{(w_j)_{j\in\N}\in\Psi_{q\bmod n}}\enspace,\]
which is a disjoint union from Lemma~\ref{clock}, and has the advantage of being a subshift, since $\sigma(\Psi_{q\bmod n})=\Psi_{q+1\bmod n}$ and $\sigma((H(q+j\bmod2n))_{j\in\N})=(H(q+1+j\bmod2n))_{j\in\N}$. 
The decoding
\[\appl\phi\Psi{B^\N}{(w_j,H(q+j\bmod2n))_{j\in\N}}{\phi_{H(q\bmod n)}((w_j)_{j\in\N})}\]
is a factorization deriving from the $n$-block map $((w_0)_j\ldots(w_{n-1})_j,H(q+j\bmod2n))_{0\le j<n}\mapsto(w_{H(q\bmod n)})\suboi0n$, since $\phi(\Psi)\supset\phi_0(\Psi_0)=B^\N$.
By definition of $\Psi_q$ and Lemma~\ref{clock}, the projections are
in $B^\omega$.
Last point, we can see that $\Psi$ is an $n$-SFT, since:
\begin{eqnarray*}
&&((w_0)_j\ldots(w_{4n-1})_j)_{j\in\N}\in\Psi\\&\Leftrightarrow&\exists q\in\ointerv0{2n},\both{((w_0)_j\ldots(w_{n-1})_j)_{j\in\N}\in\Psi_{q\bmod n}\\((w_n)_j\ldots(w_{4n-1})_j)_{j\in\N}=(H(q+j\bmod2n))_{j\in\N}}\\
&\equi&\both{\forall p\in\ointerv0n,w_p\in\sigma^{q-p\bmod n}(B^\omega)\\\forall j\in\N,(w_n)_j\ldots(w_{4n-1})_j=H(q+j\bmod2n)}\\&&\qquad\ou q=H^{-1}((w_n)_0\ldots(w_{4n-1})_0)\\
&\equi&\forall j\in\N,\both{\forall p\in\ointerv0n,(w_p)\suboi0n\in\sigma^{q-p\bmod n}(B^\omega)\suboi0n\\(w_n)_1\ldots(w_{4n-1})_1=H(q+1\bmod2n)}\\&&\qquad\ou q=H^{-1}((w_n)_0\ldots(w_{4n-1})_0)
\end{eqnarray*}
\end{prf}

\begin{prf}[of Lemma~\ref{artificiel}]
Let $\Sigma$ be a factor subshift of a SFT $\Gamma$ on alphabet $B$
such that $B^\omega\subset\Sigma$.
By Lemma~\ref{multi}, the full shift $B^\N$, where $B\subset A^n$, is a factor of a SFT $\Psi$ on alphabet $A^{4n}$ such that $\pi(\Psi)=\orb_\sigma(B^\omega)$. Let $\phi:\Psi\to B^\N$ the corresponding factorization.
$\Gamma$ is a factor of the subshift $\Gamma'\defeq\phi^{-1}(\Gamma)$.
Of course, $\pi(\Gamma')\subset\pi(\Psi)=\orb(B^\omega)$.
Moreover, it is a SFT too, since $w\in\Gamma'\equi w\in\Psi\et\phi(w)\in\Gamma$.
To sum up, $\Sigma$ is a factor of the SFT $\Gamma'$ such that
$\pi(\Gamma')\subset\Sigma$. Hence, by Proposition~\ref{factt1}, $\Sigma$ is T1.
\end{prf}

\section{Proofs of Section~\ref{k2trace}}

\begin{remark}\label{fix0}
From the definition of $w$, we immediately notice that border words have at least one letter that is not in $\phi(A)$:
$[\Upsilon]\cap[\phi(A)^{\length w}]_{\length w}=[\Upsilon]\cap[\phi(A)^{\length w}]_{2\length w}=\emptyset$.
\end{remark}

Here is a formalization of the property that border words cannot overlap each other too much:
\begin{definition}[Freezingness]
A language $W\subset A^h$ is \emph{$k$-freezing}, for some integers $k,h\in\N$, if cylinders $[W]$ and $[W]_i$ do not intersect for any $i\in\interv1k$.
\end{definition}
\begin{remark}
$W$ is $k$-freezing if and only if $[W]_i\cap[W]_j=\emptyset$ if $\abs{j-i}\in\interv1k$ (we will not have overlapping border words sharing more than $k$ letters).
\end{remark}

\begin{lemma}\label{fix1}
The set $\Upsilon$ of borders is {$(k+3\length w)$-freezing}.
\end{lemma}
\begin{prf}
~\begin{itemize}
\item If $\length w\le i\le k+3\length w$, Remark~\ref{fix0}
  gives us $[\Upsilon]\cap[\Upsilon]_i\subset[\phi(A)^{k+3\length
  w}]_{3\length w}\cap[\Upsilon]_i\subset[\phi(A)^{\length
  w}]_{i+\length w}\cap[\Upsilon]_i=\emptyset$.
\item Suppose there are an integer $i\in\ointerv0{\length w}$, words $u,v\in\orb_\gamma(w)$ and a configuration $x\in[\phi(A)^{\length w}u\overline u\phi(A)^{k+\length w}]\cap[\phi(A)^{\length w}v\overline v\phi(A)^{k+\length w}]_i$. Here we use the symmetry of words $u\overline u$ and $v\overline v$.
Let $p\defeq\min\set j{\ointerv0l}{x_{j+\length v}\notin\phi(A)}$. On the one hand:
\begin{equation}\label{eqmin}\begin{array}{lll}
p
&=&\min\set j{\ointerv0{\length w}}{u_j\notin\phi(A)}\\
&=&\length w-\max\set j{\ointerv0{\length w}}{\overline u_j\notin\phi(A)}\\
&=&2\length w-\max\set j{\ointerv0l}{x_j\notin\phi(A)}.
\end{array}\end{equation}
On the other hand:
\begin{equation}\label{eqmax}\begin{array}{lll}
p
&=&i+\min\set j{\ointerv0{\length w}}{w_j\notin\phi(A)}\\
&=&i+\length w-\max\set j{\ointerv0{\length w}}{\overline w_j\notin\phi(A)}\\
&=&i+2\length w-\max\set j{\ointerv0l}{x_j\notin\phi(A)}.
\end{array}\end{equation}
Combining Equalities \eqref{eqmin} and \eqref{eqmax}, we get $i=0$.
\end{itemize}\end{prf}

The following lemma grants the columns that it produces are in $\Sigma$.
\begin{lemma}\label{borderss}
For every border $b\in\Upsilon$, and every column $i\in\ointerv0l$, the infinite word $((\Delta_\Upsilon^j(b))_i)_{j\in\N}$ is in subshift $\Sigma$.
\end{lemma}
\begin{prf}
Consider a border word $b=a^{\length w}\gamma^j(w)\overline{\gamma^j(w)}a^{k+3\length w}$.
\begin{itemize}
\item If $i\in\ointerv0{\length w}\cup\ointerv{3\length w}l$, then $((\Delta_\Upsilon^j(b))_i)_{j\in\N}=(\phi^j(a))_{j\in\N}\in\Sigma$.
\item If $i\in\ointerv{\length w}{3\length w}$, then by a direct
  induction one finds that:
\[((\Delta_\Upsilon^j(b))_i){j\in\N}=\sigma^{\length w-\big|2\length w-i+j\big|}(w^\omega)\in\Sigma\enspace.\]
\end{itemize}\end{prf}

\begin{remark}\label{fixfix}
$\Theta$ is $(h-1)$-freezing.
\end{remark}

The following lemma grants that in particular we will be able to apply function $\Delta$ if there are two macrocells in the neighborhood.
\begin{lemma}\label{ddd}
$B\Upsilon B\Upsilon\subset\Theta$.
\end{lemma}
\begin{prf}
First, if $0<i\le k+3\length w$, then $[B\Upsilon B\Upsilon\cap A^iB\Upsilon A^{h-i}]\subset[\Upsilon]_k\cap[\Upsilon]_{k+i}=\emptyset$ (from Lemma~\ref{fix1}).

Similarly, if $k+3\length w\le i<h$, then $[B\Upsilon B\Upsilon\cap A^iB\Upsilon A^{h-i}]\subset[\Upsilon]_{h+k}\cap[\Upsilon]_{i+k}=\emptyset$.
\end{prf}

\begin{prf}[of Lemma~\ref{stabl}]
From the definition of $f$ in execution mode and of $\Delta$, we notice that:
\begin{equation}
F([\Theta])\subset[B\Upsilon]\label{stabl1}
\end{equation}
The apparition of a letter which is not in $\phi(A)$ would mean we are in the left part of an execution mode around that cell: if $x\in\az$ is such that $F(x)_0\notin\phi(A)$, then there is a $j\in\ointerv0{k+3\length w}$ such that 
$x\in[\Theta]_{-j}$.
\begin{equation}\label{stabl0}
  F^{-1}(\compl{[\phi(A)]})\subset\bigcup_{0\le j<k+3\length w}[\Theta]_{-j}
\end{equation}
Now if a configuration $x$ has an image in cylinder $[B\Upsilon]$, then, from Remark~\ref{fix0}, there is a cell $g\in\ointerv{k+\length w}{k+2\length w}$ such that $F(x)_g\notin\phi(A)$. Hence, combining with Equation~\eqref{stabl0}, with $i\defeq g-j$:
\begin{equation}\label{stabl2}
  F^{-1}([B\Upsilon])\subset\bigcup_{-2\length w<i<k+2\length w}[\Theta]_i
\end{equation}
Let $x\in F^{-1}([B\Upsilon])$. Then $x\in[\Theta]_i$ for some $i\in\oointerv{-2\length w}{k+2\length w}$ (from Equation~\eqref{stabl2}).
Then $F(x)\in[B\Upsilon]\cap F([\Theta]_i)\subset[\Upsilon]_k\cap[\Upsilon]_{k+i}$
(from Equation~\eqref{stabl1}). $\Upsilon$ being $(k+3\length w)$-freezing (Lemma~\ref{fix1}), we can conclude that $i=0$.
\begin{equation}\label{stabl3}
  F^{-1}([B\Upsilon])\subset[\Theta]
\end{equation}
Combining with Equation~\eqref{stabl1}, we get:
\begin{equation}
  F^{-1}([B\Upsilon])=[\Theta]
\end{equation}
Another consequence of Equation~\eqref{stabl3} is the stability:
\begin{equation}\label{stabl5}
  F(\compl{[\Theta]})\subset\compl{[B\Upsilon]}\subset\compl{[\Theta]}
\end{equation}
Let $x\in[\Theta]$. In particular $x\notin[\Theta]_i$ for $0<i<h$ 
Hence, $F(x)\notin[B\Upsilon]_i$ (from Equation~\eqref{stabl3}), and from Equation~\eqref{stabl1}, $F(x)\in[B\Upsilon]$. Finally, $F(x)\in[\Theta]$.
\begin{equation}\label{stabl4}
  F([\Theta])\subset[\Theta]
\end{equation}
\end{prf}

\begin{lemma}\label{inclsft2}
Let $y\in B^\Z$, and $x\in\az$ such that $\forall i\in\Z,x\suboi{ih}{ih+k}=y_i\in B$ and $x\suboi{ih+k}{(i+1)h}\in\Upsilon$. Then for $0\le q<k$, $T^q_F(x)=\pi_q(\tau_G(y))$.
\end{lemma}
\begin{prf}
We can prove by induction on the generation $j\in\N$, that for any $i\in\N$, $F^j(x)\in[G^j(y)\Upsilon]_{ih}$.
This property holds for $j=0$.
Now suppose it is true at generation $j\in\N$, and let us prove it for time $j+1$. Let $i\in\N$. $F^j(x)\in[B\Upsilon B\Upsilon]_{ih}\subset\Theta$ from the induction hypothesis and Lemma~\ref{ddd}. Therefore we are in execution mode between cells $ih$ and $(i+1)h$:
\eqn{F^{j+1}(x)\suboi0h&=&\Delta(F^j(x)\suboi{ih}{(i+1)h},F^j(x)\suboi{(i+1)h}{(i+2)h})\\
&=&\Delta(F^j(y)_i,F^j(y)_{i+1})\\
&=&F^{j+1}(y)_i\enspace.}
In particular, $T_F(x)\defeq(F^j(y)_0)_{j\in\N}=\pi_0(\tau_F(y))$.
\end{prf}

\begin{lemma}\label{inclsft1}
The trace of every configuration $x\in\az$ is in $\Sigma$.
\end{lemma}
\begin{prf}~\begin{itemize}
\item First, if $x\notin[\Theta]_{-i}$ for any $i\in\ointerv0h$, then we will always (see Lemma~\ref{stabl}) apply the default mode in cell $0$: $\forall i\in\ointerv0h,\forall j\in\N,F^j(x)\notin[\Theta]_{-i}$, and by a trivial recurrence, $T_F(x)\defeq(F^j(x)_0)_{j\in\N}=(\phi^j(x_0))_{j\in\N}\in\Sigma$.
\item If $x\in[\Theta]_{-i}$ for some $i\in\ointerv kh$, then from Lemma~\ref{inclsft2}, $\forall j\in\N,F^j(x)\in[\Theta]_{-i}\subset[\Upsilon]_{k-i}$, and by recurrence, $F^j(x)_0=\Delta_\Upsilon^j(x\suboi{k-i}{h-i})_{i-m}$, so $T_F(x)\defeq(F^j(x)_0)_{j\in\N}=(\Delta_\Upsilon^j(x\suboi{k-i}{h-i})_{i-m})_{j\in\N}\in\Sigma$ by Lemma~\ref{borderss}.
\item If $x\in[\Theta]_{qh-i}$ for some $i\in\ointerv0k$ and every $q\in\N$, then from Lemma~\ref{inclsft2}, $T_F(x)\in\stackrel\circ\tau(G)$.
\item Otherwise, there is some $i\in\ointerv0k$ and some $q\in\N^*$ such that $x\in[\Theta^q\compl\Theta]_{-i}$, and thanks to Lemma~\ref{stabl}, $\forall j\in\N,F^j(x)\in[\Theta^q\compl\Theta]_{-i}$. Let $y\in B^\N$ a configuration such that for $0\le p<q$, $x\suboi{ph-i}{ph+k-i}=y_p$ and for $p\ge q$, $y_py_{q-1}$. Then we can show by induction on $j\in\Z$ that for $0\le p<q$, $F^j(x)\suboi{ph-i}{ph+k-i}=G^k(y)_i$, from the definition of $f$. Hence, $T_F(x)=\pi_i(T_G(y))$.
\end{itemize}\end{prf}

\section{Proofs of examples}

\begin{prf}[of Example~\ref{ctrex}]
By contradiction, assume that such a CA $F$ exists. Let $f$ be the
corresponding local rule, $m$ its anchor and $d$ its diameter. Being surjective (we can immediately check that for any configuration $x\in\az$, $F^2(x)=x$), it is balanced \ie $\forall a\in A,\card{f^{-1}(a)}=\card A^{d-1}$. In particular, $\card{f^{-1}(1)}=\card{A^m0A^{d-m-1}}$. Moreover, from the definition of $\Sigma$, $f^{-1}(1)\subset A^m0A^{d-m-1}$. Equality of cardinals gives $A^m0A^{d-m-1}=f^{-1}(1)$. Hence $0^\omega\notin\Sigma$, which is a contradiction.
\end{prf}

\begin{prf}[of Example~\ref{T1T3nonT2}]
$10^*$ and $0^*1$ are included in the language of $(0^*1+1^*)0^\omega$, but not $10^*1$; hence $\Sigma$ is not a SFT.
$\Sigma$ contains two distinct transitive subshifts, namely $0^\omega$ and $1^\omega$, both of which are finite; hence it is not T2.
If $\phi$ is defined by $\forall a\in\{0,1\},\phi(a)=0$, then $0^\omega,1^\omega\in\Sigma$ and $w\defeq1\in\{0,1\}^*\setminus0^\omega$. Hence $\Sigma$ is T3.
Finally, build a 2-SFT $\Gamma\defeq((0,1)^*(1,1)+(1,1)^*)(0,0)^\omega\subset 
(\{0,1\}^2)^\omega$, and remark that $\pi_0(\Gamma)=\Sigma$; therefore
$\Sigma$ is T1.
\end{prf}

In order to prove Example~\ref{ctrex2}, we first need to prove the following lemma.
\begin{lemma}\label{clmctrx}
Let $A$ and $F$ be as in Example~\ref{ctrex2}. Then
\[
\lang(\tau(F))\cap lb^*rb^*lb^*r=\displaystyle\bigcup_{p,q\in\N}\left(lb^{2p}rb^{2q}lb^{2p}r\right).
\]
\end{lemma}
\begin{prf}
Remark that the the wall $w$ is invariant and blocking. Now we can prove that the set of configurations whose trace has a prefix in
$lb^*rb^*lb^*r$ is
\[
T_F^{-1}(lb^*rb^*lb^*r A^\omega)=\displaystyle\bigcup_{p,q\in\N}[wb^plb^qw]_{-p-1}.
\]
Intuitively, if there were no wall on the left, no particle would come from the left, since right particles cannot cross left particles. Similarly, 
there must be a wall on the right.\\
Fix $p,q\in\N$. For any $j\in\N$, 
let $i\defeq j-p-1\bmod2(p+q+1)$ and $S_j=F^j([wb^plb^qw]_{-p-1})$.
By induction on $j$ one can easily prove that
$S_j\subseteq [wb^ilb^{p+q-i}w]_{-p-1}$ if $i\le p+q$ and 
$S_j\subseteq[wb^{2(p+q)+1-i}lb^{i-p-q-1}w]_{-p-1}$ otherwise. 
Hence, the trace of a configuration $x\in[wb^plb^qw]_{-p-1}$ is $T_F(x)=\left(lb^{2p}rb^{2q}\right)^\omega$.
\end{prf}
\begin{prf}[of Example~\ref{ctrex2}]
By Lemma~\ref{clmctrx} we have that
\[
L=\lang(\tau(F))\cap lb^*rb^*lb^*r=
\displaystyle\bigcup_{p,q\in\N}\left(lb^{2p}rb^{2q}lb^{2p}r\right).
\]
Since $lb^*rb^*lb^*r$ is regular and $L$ is not, we conclude 
that $\lang(\tau(F))$ is irregular.
\end{prf}
\end{document}